\documentclass[12pt,leqno]{article}

\def\re{\mathop{\rm Re}}
\def\Vol{\mathop{\rm Vol}}

\newtheorem{theorem}{Theorem}
\newtheorem{lemma}[theorem]{Lemma}
\newtheorem{proposition}[theorem]{Proposition}
\newtheorem{sublemma}[theorem]{Sublemma}
\newtheorem{definition}[theorem]{Definition}
\newtheorem{corollary}[theorem]{Corollary}
\newtheorem{problem}[theorem]{Problem}
\newtheorem{remark}[theorem]{Remark}
\newtheorem{claim}[theorem]{Claim}
\newtheorem{assumptions}[theorem]{Assumptions}
\newtheorem{examples}[theorem]{Examples}
\newtheorem{basicfact}[theorem]{Basic Fact}

\newcommand{\begintheorem}{\addtocounter{equation}{1}\begin{theorem}}
\newcommand{\beginlemma}{\addtocounter{equation}{1}\begin{lemma}}
\newcommand{\beginproposition}{\addtocounter{equation}{1}\begin{proposition}}
\newcommand{\beginsublemma}{\addtocounter{equation}{1}\begin{sublemma}}
\newcommand{\begindefinition}{\addtocounter{equation}{1}\begin{definition}}
\newcommand{\begincorollary}{\addtocounter{equation}{1}\begin{corollary}}
\newcommand{\beginproblem}{\addtocounter{equation}{1}\begin{problem}}
\newcommand{\beginremark}{\addtocounter{equation}{1}\begin{remark}}
\newcommand{\beginclaim}{\addtocounter{equation}{1}\begin{claim}}
\newcommand{\beginassumptions}{\addtocounter{equation}{1}\begin{assumptions}}
\newcommand{\beginexamples}{\addtocounter{equation}{1}\begin{examples}}
\newcommand{\beginbasicfact}{\addtocounter{equation}{1}\begin{basicfact}}

\begin{document}

\title{Some topics concerning integrals of derivatives and difference
quotients on metric spaces}


\author{Stephen Semmes}

\date{}

\maketitle




%



\subsubsection*{Some classical inequalities}

	Suppose that $f(x)$ is a function on ${\bf R}^n$ which is,
say, continuously differentiable, and which is equal to $0$ on the
complement of a bounded set, or at least has some kind of reasonable
decay.  A well-known formula states that
\begin{equation}
\label{formula for f(x) in terms of integral of nabla f}
	f(x) = \frac{1}{\sigma_{n-1}} \sum_{j=1}^n 
			\int_{{\bf R}^n} \frac{\partial f}{\partial y_j}(y)
				\, \frac{x_j - y_j}{|x-y|^n} \, dy,
\end{equation}
where $\sigma_{n-1}$ denotes the surface measure of the unit sphere
${\bf S}^{n-1} = \{z \in {\bf R}^n : |z| = 1\}$, and $|z|$ denotes the
standard Euclidean norm of $z \in {\bf R}^n$.  This formula is given
in equation (18) on p125 of \cite{Stein1}, and it is proved by using
the Fundamental Theorem of Calculus to first write $f(x)$ as the
integral along any ray emanating from $x$ of the directional
derivative of $f$ in the direction of the ray, and then averaging over
all such rays.

	In particular, 
\begin{equation}
\label{inequality for f(x) in terms of integral of |nabla f|}
	|f(x)| \le \frac{1}{\sigma_{n-1}} \int_{{\bf R}^n} |\nabla f(y)|
					\, \frac{1}{|x-y|^{n-1}} \, dy,
\end{equation}
where $\nabla f$ denotes the gradient of $f$.  There are numerous
relatives of this inequality, such as Sobolev and Poincar\'e
inequalities, concerning the behavior of a function in terms of
integrals of its gradient, and isoperimetric-type estimates, in which
the volume of a region (or its complement) is bounded in terms of the
surface measure of its boundary.

	For instance, if $f$ is a continuously-differentiable function
on ${\bf S}^n$, then one can show that
\begin{equation}
\label{inequality on S^n, nabla f}
	\int_{{\bf S}^n} |f(x) - {\rm Av}(f)| \, dx
		\le c(n) \int_{{\bf S}^n} |\nabla f(y)| \, dy,
\end{equation}
where $c(n)$ is a positive constant that depends only on the
dimension, the $dx$ and $dy$ in the integrals denote standard surface
measure, and ${\rm Av}(f)$ is the average of $f$ over ${\bf S}^n$.  It
is understood that the gradient $\nabla f$ of $f$ refers to the
gradient in the tangential directions to ${\bf S}^n$.  Like the
previous inequality, this one can be approached in terms of the
Fundamental Theorem of Calculus and suitable averages of the estimates
that it yields.

	One can consider similar matters on spaces which are in some
way ``approximately Euclidean'' but may not be quite standard.  Now
we shall turn our attention to situations which are fractal.

\subsubsection*{A family of (fractal) Carnot--Carath\'eodory geometries}

	Suppose that $n$ is an odd integer, $n = 2m + 1$, and that $m
\ge 1$.  Let us think of ${\bf S}^n$ as being the unit sphere in ${\bf
C}^{m+1}$, i.e., ${\bf S}^n = \{w \in {\bf C}^{m+1} : |w| = 1\}$.  At
each point $z$ in ${\bf S}^n$, one has the usual real tangent space to
${\bf S}^n$, consisting of vectors pointing in directions tangent to
${\bf S}^n$ at $z$.  This can be described explicitly as
\begin{equation}
\label{real tangent space to S^n at z}
	\Bigl\{v \in {\bf C}^{m+1} : 
		\re \sum_{j=1}^{m+1} v_j \, \overline{z_j} = 0\Bigr\},
\end{equation}
where $\re a$ and $\overline{a}$ denote the real part and complex
conjugate of a complex number $a$, and $v_j$, $z_j$ denote the $j$th
components of $v, z \in {\bf C}^{m+1}$.  This is a real linear subspace
of ${\bf C}^{m+1}$ of dimension $2m + 1 = n$, which is the real dimension
of ${\bf C}^{m+1}$ minus $1$.  Inside this subspace is a complex
linear subspace of complex dimension $m$, which is given by
\begin{equation}
\label{complex subspace of real tangent space to S^n at z}
	\Bigl\{v \in {\bf C}^{m+1} : 
		\sum_{j=1}^{m+1} v_j \, \overline{z_j} = 0\Bigr\}.
\end{equation}
In other words, the original tangent space consists of the elements of
${\bf C}^{m+1}$ which are orthogonal to $z$ with respect to the
standard real inner product, and this complex subspace consists of the
elements of ${\bf C}^{m+1}$ which are orthogonal to $z$ with respect
to the standard complex Hermitian inner product.  One can account for
the difference in terms of $i z$, which lies in the first and not in
the second, and which is in fact orthogonal to the second subspace.

	If $f$ is a continuously differentiable function on ${\bf
S}^n$, then let us write $\widetilde{\nabla} f(z)$ for the version of
the gradient of $f$ which includes directional derivatives of $f$ at
$z$ in ${\bf S}^n$ only in directions in the complex subspace
(\ref{complex subspace of real tangent space to S^n at z}) of the
tangent space.  It turns out that we again have that
\begin{equation}
\label{inequality on S^n, widetilde{nabla} f}
	\int_{{\bf S}^n} |f(x) - {\rm Av}(f)| \, dx
   \le \widetilde{c}(n) \int_{{\bf S}^n} |\widetilde{\nabla} f(y)| \, dy,
\end{equation}
where $\widetilde{c}(n)$ is a positive constant that depends only on
$n$.  There is a related notion of $\widetilde{\nabla} f$ on ${\bf
R}^n$, which is connected to this version on ${\bf S}^n$ through a
``Cayley transform'', and which can also be described in terms of the
$n$-dimensional Heisenberg group.  For this there are natural
counterparts of (\ref{formula for f(x) in terms of integral of nabla
f}) and (\ref{inequality for f(x) in terms of integral of |nabla f|}),
as well as other inequalities along these lines.  See \cite{Stein2},
for instance, and also Remark \ref{affine inequalities and Heisenberg,
nilpotent groups} below.

	Ordinarily we can measure distances between points in ${\bf
S}^n$ by simply taking their distance in the Euclidean metric in ${\bf
R}^{n+1}$, or, more intrinsically, by minimizing the length of the
shortest curve in ${\bf S}^n$ between the two points.  The first
measurement of distance is less than or equal to the second one, and
the second measurement of distance is bounded by a constant factor
times the first (namely, $\pi$).  In connection with the complex
subspaces of the tangent spaces and the associated version of the
gradient $\widetilde{\nabla} f$, let us define another distance
function on ${\bf S}^n$, where the distance between two points $p$,
$q$ is the infimum of the lengths of the paths between $p$ and $q$
which have the additional feature that their tangent vectors lie in
the complex subspace of the tangent space to ${\bf S}^n$ at the given
point (as in (\ref{complex subspace of real tangent space to S^n at
z})).  More precisely, in this definition let us consider curves
$\gamma(t)$ in ${\bf S}^n$, where $t$ runs through a closed interval
in ${\bf R}$, which are piecewise-smooth, and which satisfy this
condition on the tangent vectors except possibly at finitely many
values of $t$.

	It is a classical fact that for each pair of points $p$, $q$
in ${\bf S}^n$ there is such a path $\gamma(t)$.  This reflects the
``nonintegrability'' of the complex subspaces of the tangent spaces.
Thus the distance between $p$ and $q$ defined above is always finite,
and of course it is greater than or equal to the distance defined by
minimizing the lengths of all curves between $p$ and $q$, and not just
the ones with this condition on the tangent vectors.

	This new distance function defines the same topology on ${\bf
S}^n$ as the usual one, but it is quite different geometrically.
Basically, it looks like the square root of the other distances in
some directions, to balance the restriction on the curves.  Indeed, if
one starts at a point $z$ in ${\bf S}^n$ and tries to go a small
distance in the direction $i z$, which lies in the tangent space
(\ref{real tangent space to S^n at z}) to ${\bf S}^n$ at $z$ but not
in the complex subspace (\ref{complex subspace of real tangent space
to S^n at z}), then one cannot use curves which go directly in this
direction.  Instead one uses more indirect paths which follow the
complex subspaces of the tangent spaces.

	The new distance function in fact defines a geometry on ${\bf
S}^n$ which is fractal.  For balls of radius $r \le 1$, the measure of
the ball is on the order of $r^{n+1}$, rather than $r^n$, as for the
original geometry.  In the new geometry, ${\bf S}^n$ has Hausdorff
dimension $n+1$, rather than $n$.

	These changes affect the kind of inequalities for functions
$f$ and the modified gradient $\widetilde{\nabla} f$ that one gets.
This does not come up so much in (\ref{inequality on S^n, nabla f}),
(\ref{inequality on S^n, widetilde{nabla} f}), except for the
constant, but the distance function and Hausdorff dimension are
explicitly involved in other inequalities.

\beginremark
\label{affine inequalities and Heisenberg, nilpotent groups}
{\rm In \cite{LYZ} ``affine'' versions of Sobolev inequalities are
discussed, i.e., inequalities which only involve the vector space
structure and Lebesgue measure on ${\bf R}^n$, and not an inner
product or norm.  To my knowledge analogous matters have not been
investigated for the Heisenberg groups or other (noncommutative
simply-connected) nilpotent Lie groups.  That is, one would be
interested in inequalities that use the underlying Haar measure (which
is Lebesgue measure in standard coordinates), group structure, and
dilations, but not otherwise a choice of norm or metric.  To look at
it another way, in place of invariance under invertible linear
transformations on ${\bf R}^n$ that preserve volume, one could try to
get invariance under automorphisms of a nilpotent Lie group which
preserve volume.  On the one hand, noncommutativity makes the group
more tricky, and on the other hand, the group of automorphisms is smaller.
}
\end{remark}

\subsubsection*{Other geometries}

	What about other kinds of geometry?  As a basic set-up, let
$(M, d(x,y))$ be a metric space, and let $\mu$ be a reasonable measure
on $M$.  If $f$ is a function on $M$, define $D_\epsilon f$ for
$\epsilon > 0$ by
\begin{equation}
\label{def of D_epsilon f}
	D_\epsilon f(x) = 
   \sup \biggl\{\frac{|f(y) - f(x)|}{\epsilon} : 
				d(x,y) \le \epsilon \biggr\}.
\end{equation}
We would like to know about inequalities concerning the behavior of
$f$ and integrals involving $D_\epsilon f$.

	For example, suppose that $M$ is a fractal set in a Euclidean
space, like the Sierpinski gasket or carpet, with the usual Euclidean
metric restricted to the set as our metric.  For fractals like these
there are natural measures $\mu$ so that the whole set has positive
finite measure, and so that the measure assigns equal values to
different ``pieces'' of the set that are essentially the same in the
construction, such as intersections of the set with certain triangles
or squares from the construction.  Thus there are well-behaved ingredients
for the basic set-up.

	However, the inequalities do not work in the same way as
before.  One can have functions $f$ on the Sierpinski gasket or
carpet, for instance, which are equal to $0$ on one fixed (nontrivial)
part of the set, are equal to $1$ on another fixed part of the set,
and for which the integral of $D_\epsilon f$ tends to $0$ as
$\epsilon$ tends to $0$.  This happens by concentrating the
oscillations of $f$ in sufficiently-small sets, which is possible for
these fractal spaces.  For these types of spaces adjustments are
needed in the inequalities to be considered in order to take the
geometry into account.  See \cite{Kigami}.

	In \cite{BP1, Laakso1}, remarkable families of metric spaces
are discussed in which one does have inequalities as before, i.e.,
with suitable bounds for a function $f$ in terms of integrals of
$D_\epsilon f$ which remain uniform as $\epsilon \to 0$.

	It may seem surprising that there are genuinely fractal metric
spaces with approximately the same kind of inequalities concerning the
behavior of a function $f$ and integrals of $D_\epsilon f$ (with
uniform bounds in $\epsilon$ in a suitable range).  However, there are
basic differences in the ``internal'' structure of the spaces,
regarding intermediate dimensions (compared to the dimension of the
space).  See \cite{misha-cc} in this connection, for the fractal
geometries on the spheres ${\bf S}^n$ described earlier, and other
analogous geometries.

	There are related issues in non-fractal contexts, such as
weighted inequalities on ${\bf R}^n$.  For (nonsmooth) conformal
deformations associated to ``strong-$A_\infty$ weights'', as in
\cite{DS-sob}, there is good behavior internally, and a form of this
is described in \cite{Se-barc}.  On the other hand, there are weights
for which there are plenty of nice inequalities for functions and
their derivatives, but which do not behave so well internally.  For
instance, the weight might be a continuous function on ${\bf R}^3$
which vanishes on a circle and is positive otherwise.

\subsubsection*{Graphs}

	Suppose that we have a graph, which is to say a nonempty set
$V$ of vertices and a set $E$ of edges, which can be represented by
unordered pairs of distinct vertices in $V$.  It is natural to that
ask the graph has locally bounded geometry, in the sense that there is
a constant $C_0$ so that each element of $V$ is adjacent to at most
$C_0$ elements of $V$.  If $C_0$ is not too large, then this is also a
significant condition when $V$ is finite.

	In addition we ask that the graph is connected, which means
that for any pair of points $x$, $y$ in $V$ there is a path in $V$
which goes from $x$ to $y$, i.e., a finite sequence of successively
adjacent elements of $V$ connecting $x$ to $y$.  Thus we can define
the distance $d(x,y)$ between $x$ and $y$ to be the length of the
shortest such path, where the length of the path is the number of
elements of $V$ in the sequence minus $1$.

	If $x$ lies in $V$, let us write $N(x)$ for the set of elements
of $V$ which are adjacent to $x$.  For a real or complex-valued function
$f$ on $V$, define $D(f)$ on $V$ by
\begin{equation}
	D(f)(x) = \max \{|f(w) - f(x)| : w \in N(x)\}.
\end{equation}
This plays the role of the absolute value of the gradient in the
present discussion, and is the same as $D_1(f)$ in (\ref{def of
D_epsilon f}).  In place of integrals before, one can now use sums.

	For example, for $V$ one can take ${\bf Z}^n$, the set of
$n$-tuples of integers, and one can say that two elements of ${\bf
Z}^n$ are adjacent when there difference has one component equal to
$\pm 1$ and the others equal to $0$.  One might instead take $V$ to be
the set of $n$-tuples of integers which lie in a box, with the same
condition for adjacency.

	Another class of examples comes from the discrete version of
the Heisenberg groups.  Fix a positive integer $n$, and let $H_n$ be
the group with generators $a_1, \ldots, a_n$, $b_1, \ldots, b_n$,
subject to the relations that the $a_i$'s commute with each other, the
$b_i$'s commute with each other, the $a_i$'s commute with the $b_j$'s
when $i \ne j$, and the commutators $b_i a_i b_i^{-1} a_i^{-1}$ are
equal to each other and commute with all generators of the group.  Let
us write $c$ for the common value of the commutators $b_i a_i b_i^{-1}
a_i^{-1}$.

	There is a standard normal form for the elements of $H_n$, 
namely
\begin{equation}
\label{normal form in H_n}
	a_1^{r_1} \cdots a_n^{r_n} b_1^{s_1} \cdots b_n^{s_n} c^t,
\end{equation}
where $r_1, \ldots, r_n$, $s_1, \ldots, s_n$, and $t$ are integers.
That is, every element of $H_n$ can be written in this form in a
unique manner.

	Let us say that two elements $x$, $y$ of $H_n$ are adjacent
if $y$ can be written as $x u$, where $u$ is either a generator or its
inverse.  This defines the \emph{Cayley graph} associated to $H_n$,
and we obtain from it a distance function $d(x,y)$ on $H_n$.  This
distance function is invariant under left translations on $H_n$.

	One can show that the number of elements of $H_n$ whose
distance to the identity element is less than or equal to $\ell$ is of
the order $\ell^{2n+2}$ when $\ell$ is a positive integer (rather than
$\ell^{2n+1}$, as might be suggested by the normal form (\ref{normal
form in H_n}), and as would be the case for ${\bf Z}^{2n+1}$).  In
this regard a basic ingredient is the identity
\begin{equation}
	b_i^j a_i^j b_i^{-j} a_i^{-j} = c^{j^2}.
\end{equation}

	A very nice approach to various inequalities for $H_n$ and other
groups is given in \cite{CSC}.

\subsubsection*{Isoperimetric inequalities}

	The classical isoperimetric inequality states that if $\Omega$
is a region in ${\bf R}^n$ with finite volume, then
\begin{equation}
\label{Vol_n(Omega) le c_n Vol_{n-1}(partial Omega)^{n/(n-1)}}
	{\Vol}_n(\Omega) \le c_n \, {\Vol}_{n-1}(\partial \Omega)^{n/(n-1)},
\end{equation}
where $c_n$ is a positive constant which is chosen so that equality
holds when $\Omega$ is a ball.  Here $\partial \Omega$ denotes the
boundary of $\Omega$, and $\Vol_n$, $\Vol_{n-1}$ denote the $n$ and
$(n-1)$-dimensional volumes (of a given set).  One might as well assume
that $\Omega$ is a bounded open set with reasonably-nice boundary,
so that the volumes in question can be defined as in vector calculus,
but there are more elaborate versions using Hausdorff measures or
generalized derivatives and perimeter which apply more generally.

	Let us note that the power $n/(n-1)$ on the right side of
(\ref{Vol_n(Omega) le c_n Vol_{n-1}(partial Omega)^{n/(n-1)}}) is
determined by considerations of scaling.  That is, if we apply a
dilation $x \mapsto \lambda \, x$ on ${\bf R}^n$, where $\lambda$ is a
positive real number, then the $n$-dimensional volume of a set after
dilation is equal to $\lambda^n$ times the $n$-dimensional volume of
the set before dilation, and similarly the $(n-1)$-dimensional volume
of a set after dilation is equal to $\lambda^{n-1}$ times the
$(n-1)$-dimensional volume of the set before dilation.  The power
$n/(n-1)$ on the right side of (\ref{Vol_n(Omega) le c_n
Vol_{n-1}(partial Omega)^{n/(n-1)}}) ensures that we get the same
factor of $\lambda^n$ on both sides of (\ref{Vol_n(Omega) le c_n
Vol_{n-1}(partial Omega)^{n/(n-1)}}) after dilation.

	If we do not ask for the sharp constant, then
(\ref{Vol_n(Omega) le c_n Vol_{n-1}(partial Omega)^{n/(n-1)}}) can be
obtained from (\ref{inequality for f(x) in terms of integral of |nabla
f|}), as follows.  Let $f(x)$ be the characteristic function of
$\Omega$, which is to say that $f(x) = 1$ when $x$ lies in $\Omega$
and $f(x) = 0$ when $x$ lies in the complement of $\Omega$ in ${\bf
R}^n$.  Of course $f$ is not a continuous function, let alone smooth,
but its ``generalized gradient'' is given by the unit normal vector on
$\partial \Omega$ pointing into $\Omega$ times surface measure on
$\partial \Omega$.  The appropriate version of (\ref{inequality for
f(x) in terms of integral of |nabla f|}) becomes
\begin{equation}
\label{|f(x)| le (potential along partial Omega)}
	|f(x)| \le \frac{1}{\sigma_{n-1}} \int_{\partial \Omega} 
				\frac{1}{|x-y|^{n-1}} \, d{\Vol}_{n-1}(y).
\end{equation}
This inequality can be derived through the same kind of computation as
before, or it can be viewed as a consequence of the previous one,
using suitable approximations of $f$ by smooth functions.

	If $\mu$ is any nonnegative measure on ${\bf R}^n$ with finite
total mass, then there is a ``weak-type'' inequality to the effect
that for each $t > 0$ the ordinary volume of the set
\begin{equation}
\label{set of x in R^n where potential is ge t}
   \biggl\{x \in {\bf R}^n : \frac{1}{\sigma_{n-1}} \int_{{\bf R}^n}
			\frac{1}{|x-y|^{n-1}} \, d\mu(y) \ge t \biggr\}
\end{equation}
is bounded by the product of $t^{-n/(n-1)}$, the total mass of $\mu$
raised to the power $n/(n-1)$, and a constant that depends only on
$n$.  This follows from the argument on p120 of \cite{Stein1}.
Alternatively, the normability of the ``weak-type'' space $L^{n/(n-1),
\infty}$ can be used, as in Theorem 3.21 on p204 of \cite{SW}.  In our
situation, where $f$ is the characteristic function of $\Omega$ and
satisfies (\ref{|f(x)| le (potential along partial Omega)}), we can
apply the weak-type inequality for (\ref{set of x in R^n where
potential is ge t}) with $t = 1$ to get a bound for the volume of
$\Omega$ by a constant times $\Vol_{n-1}(\partial \Omega)$.

	These are standard types of arguments which are applied in
many contexts.  In particular, there are numerous inequalities for
potentials such as the weak-type inequality just mentioned.

	Let us return now to the matter of the exponent $n/(n-1)$ on
the right side of (\ref{Vol_n(Omega) le c_n Vol_{n-1}(partial
Omega)^{n/(n-1)}}).  If $\Omega$ is a ball of radius $r$, then the
volume of $\Omega$ is equal to a constant depending only on $n$ times
$r^n$, and $\Vol_{n-1}(\partial \Omega)$ is equal to another constant
that depends only on $n$ times $r^{n-1}$.  This again determines the
exponent $n/(n-1)$.

	In general, the exponent in an isoperimetric inequality (or
other factors for related inequalities) is connected to volume growth
in the space.  In some cases, including suitable conditions of
negative curvature, there are isoperimetric inequalities with the
exponent equal to $1$ (and no additional factors), and the volume
grows exponentially.  Here we are more concerned with spaces in which
the volume of balls grows at most polynomially, and for which the
exponent in an isoperimetric inequality would be larger than $1$ (or
there would be other adjustments, such as an extra factor of the
diameter of $\Omega$).

\subsubsection*{Hyperbolic groups}

	Suppose that $G$ is a countably-infinite group which is
generated by a finite set $S$.  It will be convenient to assume that
$S$ is \emph{symmetric}, by which we mean that $S$ contains the
inverses of its elements.  As in the case of the discrete versions of
the Heisenberg groups, the \emph{Cayley graph} associated to $G$
and $S$ consists of the elements of $G$ as vertices with an edge
from $x$ to $y$ in $G$ whenever there is a generator $u$ in $S$
such that $y = x \, u$.

	Every pair of vertices in the graph can be connected by a
path, and this gives rise to a metric on $G$, in which the distance
between two points is the length of the shortest path that connects
them.  This distance function is invariant under \emph{left
translations} on the group, which are transformations of the form $x
\mapsto a \, x$, where $a$ is any element of $G$.

	\emph{Hyperbolic groups} can be defined in terms of certain
\emph{negative curvature} properties.  See \cite{CP, Ghys-Harpe,
misha-hypgroups}.  An important aspect of these groups are their
``spaces at infinity'', and indeed these groups lead to many
interesting geometries this way.

	A special case concerns cocompact (also known as uniform)
lattices associated to rank-$1$ symmetric spaces.  The spaces at
infinity corresponding to these have Euclidean geometry for classical
hyperbolic spaces (with constant negative curvature), and otherwise
have Carnot--Carath\'eodory geometry.  The geometry on the unit sphere
in ${\bf C}^{m+1}$ using the complex subspaces of the ordinary tangent
spaces discussed earlier corresponds to $(m+1)$-dimensional complex
hyperbolic space, for instance.

	The fractal spaces considered in \cite{Bourdon, BP1} are also
connected to hyperbolic groups, and to \emph{buildings} instead of
symmetric spaces.  Concerning the general subject of buildings, see
\cite{Brown}.

\subsubsection*{Norms on ${\bf R}^n$}

	Fix a positive integer $n$, and suppose that $\|\cdot\|$ is a
norm on ${\bf R}^n$.  Thus $\|x\|$ is a nonnegative real number for
each $x$ in ${\bf R}^n$ which is $0$ exactly when $x = 0$, $\|t \, x\|
= |t| \, \|x\|$ for all real numbers $t$ and $x$ in ${\bf R}^n$, and
$\|x+y\| \le \|x\| + \|y\|$ for all $x$, $y$ in ${\bf R}^n$.
Associated to this norm is the distance function $d(x,y) = \|x - y\|$
on ${\bf R}^n$.

	Erwin Lutwak once explained to me how there are many
interesting questions in this setting related to the sort of topics
and questions that we have been considering here.  The distance
function just defined and the ordinary Euclidean distance are each
bounded by constant multiples of the other, so that similar
inequalities are around as in the Euclidean case, but of course the
constants do not have to be the same, nor the geometry of various
optimizing objects, such as for isoperimetry.  For that matter, there
are numerous questions about how to make basic measurements, such as
area of submanifolds, and different choices lead to different
geometric possibilities.

	A basic reference related to these matters is \cite{Thompson}.

	In addition to ${\bf R}^n$ itself, one might consider lattices
contained therein.  Of course the geometry of a norm in connection
with a lattice can be quite remarkable.  See \cite{Cassels1}.

\subsubsection*{$p$-Adic numbers and absolute values}

	Let $p$ be a prime number.  The $p$-adic absolute value
$|x|_p$ of a rational number $x$ is defined as follows.  If $x = 0$,
then $|x|_p = 0$.  If $x = (a/b) p^k$, where $a$, $b$, and $k$ are
integers, with $a, b \ne 0$ and neither $a$ nor $b$ divisible by $p$,
then $|x|_p = p^{-k}$.  Two basic properties of $|\cdot |_p$ are
\begin{equation}
	|x + y|_p \le \max(|x|_p, |y|_p), \qquad |x \, y|_p = |x|_p \, |y|_p
\end{equation}
for all $x, y \in {\bf Q}$.  These are not difficult to verify.
Notice that
\begin{equation}
	(p-1) \sum_{j=0}^k p^j = -1 + p^{k+1}
\end{equation}
for $k \ge 0$, so that $-1$ can be approximated in the $p$-adic absolute
value by positive integers.

	Just as the real numbers can be viewed as a completion of the
rational numbers with respect to the ordinary absolute value function
$|x|$, the \emph{$p$-adic numbers} can be viewed as a completion of
the rational numbers with respect to the $p$-adic absolute value
function $|x|_p$.  For simplicity, let us restrict our attention to
rational numbers here, with measurements in terms of an absolute value
function.  See \cite{Cassels2, Gouvea} for more information about
$p$-adic numbers and absolute values.

	Suppose that we have a graph again, consisting of a nonempty
set $V$ of vertices which is at most countable, and a set $E$ of edges
which can be represented by unordered pairs of distinct vertices in
$V$.  We ask that the graph has at most finitely many edges attached
to any fixed vertex, and, as before, one might also ask for a bound on
the number of such edges.

	Let $(a_{u,v})$ be a matrix of rational numbers, where $u$ and
$v$ run through the set $V$ of vertices.  We ask that $a_{u,v}$ be
zero unless either $u = v$ or there is an edge in the graph between
$u$ and $v$.  From this matrix we get an operator $A$ on functions on
$V$ with values in ${\bf Q}$, given by
\begin{equation}
	A(f)(v) = \sum_{u \in V} a_{u,v} \, f(u).
\end{equation}
Under our assumptions, the sum on the right is a finite sum for each
$v$ in $V$.

	A fundamental feature of the operator $A$ is that
\begin{equation}
\label{sup_{v in V} |A(f)(v)| le ...}
	\sup_{v \in V} |A(f)(v)| 
		\le \Bigl(\sup_{v \in V} \sum_{u \in V} |a_{u,v}| \Bigr)
			\, \Bigl(\sup_{w \in V} |f(w)| \Bigr).
\end{equation}
Here we are using the classical absolute values on ${\bf Q}$.  It is
not hardy to verify (\ref{sup_{v in V} |A(f)(v)| le ...}), using the
triangle inequality.  This estimate is sharp, in the sense that for
each $v$ in $V$ there is a function $f$ on $V$ such that
\begin{equation}
	A(f)(v) = \sum_{u \in V} |a_{u,v}|
\end{equation}
and $f(w) = \pm 1$ for all $w$ in $V$.

	Now let $p$ be a prime number again, and consider the $p$-adic
absolute value function $|\cdot |_p$.  For this we have that
\begin{equation}
\label{sup_{v in V} |A(f)(v)|_p le ...}
	\sup_{v \in V} |A(f)(v)|_p
		\le \Bigl(\sup_{v \in V} \sup_{u \in V} |a_{u,v}|_p \Bigr)
			\, \Bigl(\sup_{w \in V} |f(w)|_p \Bigr).
\end{equation}
This estimate can be verified using the ultrametric version of the
triangle inequality for the $p$-adic absolute value function.  We also
have that (\ref{sup_{v in V} |A(f)(v)|_p le ...}) is sharp, in that
for each $v$ in $V$ there is a function $f$ on $V$ such that
\begin{equation}
	|A(f)(v)|_p = \sup_{u \in V} |a_{u,v}|_p
\end{equation}
and $f(w) = 1$ for one $w$ in $V$, $f(z) = 0$ at all other points in $V$.

	Another feature of the operator $A$ is that
\begin{equation}
\label{sum_{v in V} |A(f)(v)| le ...}
	\sum_{v \in V} |A(f)(v)|
		\le \Bigl(\sup_{u \in V} \sum_{v \in V} |a_{u,v}| \Bigr)
			\, \Bigl(\sum_{w \in V} |f(w)| \Bigr)
\end{equation}
and
\begin{equation}
\label{sum_{v in V} |A(f)(v)|_p le ...}
	\sum_{v \in V} |A(f)(v)|_p
		\le \Bigl(\sup_{u \in V} \sum_{v \in V} |a_{u,v}|_p \Bigr)
			\, \Bigl(\sum_{w \in V} |f(w)|_p \Bigr)
\end{equation}
for rational-valued functions on $V$.  These inequalities are not
difficult to verify.  Notice that if $y$ is an element of $V$ and
$f_y$ is the function on $V$ defined by $f_y(y) = 1$, $f_y(z) = 0$
when $z \ne y$, then $A(f_y)(v) = a_{y,v}$,
\begin{equation}
	\sum_{v \in V} |A(f_y)(v)| = \sum_{v \in V} |a_{y,v}|, \quad
		\sum_{v \in V} |A(f_y)(v)|_p = \sum_{v \in V} |a_{y,v}|_p,
\end{equation}
and $\sum_{w \in V} |f_y(w)| = 1$, $\sum_{w \in V} |f_y(w)|_p = 1$.
This shows the sharpness of (\ref{sum_{v in V} |A(f)(v)| le ...}),
(\ref{sum_{v in V} |A(f)(v)|_p le ...}).

\subsubsection*{Other types of computations on graphs}

	There are contexts in which one looks at different kinds of
operations associated to vertices or edges in a graph.  An
instance of this is given by \emph{Boolean circuits}, with the
operations $\land$, $\lor$, and $\neg$ ({\bf and}, {\bf or}, and
negation).  See \cite{LP, P}, for example.

	Connectedness is a basic issue, i.e., is there a path from a
vertex $p$ to a vertex $q$, perhaps satisfying some auxiliary
condition.  There can also be questions of having a number of paths
from one point to another, and not just a single path.  See \cite{CS}.



\begin{thebibliography}{BakCLS}


\bibitem [Ale] {Alexopoulos} G.~Alexopoulos, {\it Sub-Laplacians with
Drift on Lie Groups of Polynomial Volume Growth}, Memoirs of the
American Mathematical Society {\bf 739}, 2002.

\bibitem [Amb] {Ambrosio} L.~Ambrosio, {\it Some fine properties of
sets of finite perimeter in Ahlfors regular metric measure spaces},
Advances in Mathematics {\bf 159} (2001), 51--67.

\bibitem [AmbM] {Ambrosio-Magnani} L.~Ambrosio and V.~Magnani, {\it Some
fine properties of BV functions on sub-Riemannian groups}, preprint, 2002.

\bibitem [AmbT] {Ambrosio-Tilli} L.~Ambrosio and P.~Tilli, {\it
Selected Topics on ``Analysis in Metric Spaces''}, Scuola Normale
Superiore, Pisa, 2000.

\bibitem [BakCLS] {BCLSC} D.~Bakry, T.~Coulhon, M.~Ledoux, and
L.~Saloff-Coste, {\it Sobolev inequalities in disguise}, Indiana
University Mathematics Journal {\bf 44} (1995), 1033--1074.

\bibitem [Bal] {Ballmann} W.~Ballmann, {\it Lectures on Spaces of
Nonpositive Curvature}, Birkh\"auser, 1995.

\bibitem [BalGS] {BalGS} W.~Ballmann, M.~Gromov, and V.~Schroeder,
{\it Manifolds of Nonpositive Curvature}, Birkh\"auser, 1985.

\bibitem [BeaGaG] {BGG} R.~Beals, B.~Gaveau, and P.~Greiner, {\it
Hamilton--Jacobi theory and the heat kernel on Heisenberg groups},
Journal de Math\'ematiques Pures et Appliqu\'ees (9) {\bf 79} (2000),
633--689.

\bibitem [BeaGr] {BG} R.~Beals and P.~Greiner, {\it Calculus on Heisenberg
Manifolds}, Annals of Mathematics Studies {\bf 119}, 1988.

\bibitem [BeaGrS] {BGS} R.~Beals, P.~Greiner, and N.~Stanton, {\it The
heat equation on a CR manifold}, Journal of Differential Geometry {\bf
20} (1984), 343--387.

\bibitem [Bel] {Bellaiche} A.~Bella\"{\i}che, {\it The tangent space
in sub-Riemannian geometry}, in {\it Sub-Riemannian Geometry},
A.~Bella\"{\i}che and J.-J.~Risler, editors, 1--78, Birkh\"auser,
1996.

\bibitem [BerV] {BV} V.~Berestovskij and A.~Vershik, {\it Manifolds
with intrinsic metric, and nonholonomic spaces}, Representation theory
and dynamical systems, Advances in Soviet Mathematics {\bf 9} (1992),
253-267.

\bibitem [Bj\"oBS] {BBS} A.~Bj\"orn, J.~Bj\"orn, and
N.~Shanmugalingam, {\it The Dirichlet problem for $p$-harmonic
functions on metric measure spaces}, Reports of the Department of
Mathematics, Link\"opings Universitet, 2001.

\bibitem [Bj\"oMS] {BMS} J.~Bj\"orn, P.~MacManus, and S.~Shanmugalingam,
{\it Fat sets and pointwise boundary estimates for $p$-harmonic functions
in metric spaces}, Journal d'Analyse Math\'ematique
{\bf 85} (2001), 339--369.

\bibitem [BonK] {BK} M.~Bonk and B.~Kleiner, {\it Quasisymmetric
parameterizations of two-dimensional metric spheres}, preprint, 2001.

\bibitem [Bou] {Bourdon} M.~Bourdon, {\it Immeubles hyperboliques,
dimension conforme et rigidit\'e de Mostow}, Geometric and Functional
Analysis {\bf 7} (1997), 245--268.

\bibitem [BouP1] {BP1} M.~Bourdon and H.~Pajot, {\it Poincar\'e
inequalities and quasiconformal structure on the boundary of some
hyperbolic buildings}, Proceedings of the American Mathematical
Society {\bf 127} (1999), 2315--2324.

\bibitem [BouP2] {BP2} M.~Bourdon and H.~Pajot, {\it Rigidity of
quasi-isometries for some hyperbolic buildings}, Commentarii
Mathematici Helvetici {\bf 75} (2000), 701-736.

\bibitem [BouP3] {BP3} M.~Bourdon and H.~Pajot, {\it Quasi-conformal
geometry and hyperbolic geometry}, in {\it Rigidity in Dynamics and
Geometry}, M.~Burger and A.~Iozzi, editors, 1--15, Springer-Verlag, 2002.

\bibitem [BouP4] {BP4} M.~Bourdon and H.~Pajot, {\it Cohomologie $\ell_p$
et espaces de Besov}, preprint, 2002.

\bibitem [Bro] {Brown} K.~Brown, {\it Buildings}, Springer-Verlag,
1989.

\bibitem [CarS] {CS} A.~Carbone and S.~Semmes, {\it A Graphic Apology
for Symmetry and Implicitness}, Oxford University Press, 2001.

\bibitem [Cas1] {Cassels1} J.~Cassels, {\it An Introduction to the
Geometry of Numbers}, Springer-Verlag, 1971.

\bibitem [Cas2] {Cassels2} J.~Cassels, {\it Local Fields}, Cambridge
University Press, 1986.

\bibitem [Che] {Cheeger} J.~Cheeger, {\it Differentiability of
Lipschitz functions on metric measure spaces}, Geometric and
Functional Analysis {\bf 9} (1999), 428--517.

\bibitem [ColM1] {ColM1} T.~Colding and W.~Minicozzi II, {\it Weyl
type bounds for harmonic functions}, Inventiones Mathematicae {\bf
131} (1998), 257--298.

\bibitem [ColM2] {ColM2} T.~Colding and W.~Minicozzi II, {\it Liouville
theorems for harmonic sections and applications}, Communications
on Pure and Applied Mathematics {\bf 51} (1998), 113-138.

\bibitem [Coo] {C} M.~Coornaert, {\it Mesures de Patterson--Sullivan
sur le bord d'un espace hyperbolique au sens de M.~Gromov}, Pacific
Journal of Mathematics {\bf 159} (1993), 241--270.

\bibitem [CooDP] {CDP} M.~Coornaert, T.~Delzant, and A.~Papadopoulos,
{\it G\'eom\'etrie et Th\'eorie des Groupes, les Groupes Hyperboliques
de M.~Gromov}, Lecture Notes in Mathematics {\bf 1441},
Springer-Verlag, 1991.

\bibitem [CooP] {CP} M.~Coornaert and A.~Papadopoulos, {\it Symbolic
Dynamics and Hyperbolic Groups}, Lecture Notes in Mathematics {\bf
1539}, Springer-Verlag, 1993.

\bibitem [Cou1] {Coulhon1} T.~Coulhon, {\it Espaces de Lipschitz et
in\'egalit\'es de Poincar\'e}, Journal of Functional Analysis {\bf 136}
(1996), 81--113.

\bibitem [Cou2] {Coulhon2} T.~Coulhon, {\it Ultracontractivity and
Nash type inequalities}, Journal of Functional Analysis {\bf 141}
(1996), 510--539.

\bibitem [Cou3] {Coulhon3} T.~Coulhon, {\it Random walks and geometry
on infinite graphs}, in {\it Lecture Notes on Analysis on Metric
Spaces}, L.~Ambrosio and F.~Serra Cassano, editors, 5--36, Scuola
Normale Superiore, Pisa, 2000.

\bibitem [CouG1] {CG1} T.~Coulhon and A.~Grigoryan, {\it On-diagonal
lower bounds for heat kernels and Markov chains}, Duke Mathematical
Journal {\bf 89} (1997), 133--199.

\bibitem [CouG2] {CG2} T.~Coulhon and A.~Grigoryan, {\it Random walks
on graphs with regular volume growth}, Geometric and Functional
Analysis {\bf 8} (1998), 656--701.

\bibitem [CouG3] {CG3} T.~Coulhon and A.~Grigoryan, {\it Pointwise
estimates for transition probabilities of random walks on infinite
graphs}, in {\it Proceedings of the Conference ``Fractals in Graz''},
to appear, 2002.

\bibitem [CouS] {CSC} T.~Coulhon and L.~Saloff-Coste, {\it
Isop\'erim\'etrie pour les groupes et les vari\'et\'es}, Revista
Matem\'atica Iberoamericana {\bf 9} (1993), 293--314.

\bibitem [DavS1] {DS-sob} G.~David and S.~Semmes, {\it
Strong-$A_\infty$ weights, Sobolev inequalities, and quasiconformal
mappings}, in {\it Analysis and Partial Differential Equations},
C.~Sadosky, editor, 101--111, Marcel Dekker, 1990.

\bibitem [DavS2] {DS-frac} G.~David and S.~Semmes, {\it Fractured
Fractals and Broken Dreams: Self-Similar Geometry through Metric
and Measure}, Oxford University Press, 1997.

\bibitem [DeVS] {DeVore-Sharpley} R.~DeVore and R.~Sharpley, {\it
Maximal Functions Measuring Smoothness}, Memoirs of the American
Mathematical Society {\bf 293}, 1984.

\bibitem [FerR] {Josechu-R} J.~Fern\'andez and J.~Rodr\'{\i}guez,
{\it Area growth and Green's function of Riemann surfaces},
Arkiv f\"or Matematik {\bf 30} (1992), 83--92.

\bibitem [FolS] {FS} G.~Folland and E.~Stein, {\it Hardy Spaces on
Homogeneous Groups}, Princeton University Press, 1982.

\bibitem [FraHK] {FHK} B.~Franchi, P.~Haj{\l}asz, and P.~Koskela, {\it
Definitions of Sobolev classes on metric spaces}, Annales de
l'Institut Fourier (Grenoble) {\bf 49} (1999), 1903--1924.

\bibitem [FraSS1] {FSSC1} B.~Franchi, R.~Serapioni, and F.~Serra
Cassano, {\it Sur les ensembles de p\'erim\`etre fini dans le groupe
de Heisenberg}, Comptes Rendus de l'Acad\'emie des Sciences
Paris S\'er.~I Math.\ {\bf 329} (1999), 183--188.

\bibitem [FraSS2] {FSSC2} B.~Franchi, R.~Serapioni, and F.~Serra
Cassano, {\it Rectifiability and perimeter in the Heisenberg group},
Mathematiche Annalen {\bf 321} (2001), 479--531.

\bibitem [FraSS3] {FSSC3} B.~Franchi, R.~Serapioni, and F.~Serra
Cassano, {\it Regular hypersurfaces, intrinsic perimeter, and implicit
function theorem in Carnot groups}, to appear, Communications in
Analysis and Geometry.

\bibitem [GhyH] {Ghys-Harpe} E.~Ghys and P.~de la Harpe, editors, {\it
Sur les Groupes Hyperboliques d'apr\`es Mikhael Gromov}, Birkh\"auser,
1990.

\bibitem [Gri1] {Grigoryan1} A.~Grigoryan, {\it Heat kernel upper
bounds on a complete non-compact manifold}, Revista Matem\'atica
Iberoamericana {\bf 10} (1994), 395--452.

\bibitem [Gri2] {Grigoryan2} A.~Grigoryan, {\it Analytic and geometric
background of recurrence and non-explosion of the Brownian motion on
Riemannian manifolds}, Bulletin of the American Mathematical Society
(N.S.) {\bf 36} (1999), 135--249.

\bibitem [Gro1] {misha-hyperbolic} M.~Gromov, {\it Hyperbolic manifolds,
groups, and actions}, in {\it Riemann Surfaces and Related Topics:
Proceedings of the 1978 Stony Brook Conference}, I.~Kra and B.~Maskit,
editors, 183--213, Annals of Mathematics Studies {\bf 97}, Princeton
University Press, 1981.

\bibitem [Gro2] {misha-green} M.~Gromov, {\it Structures M\'etriques
pour les Vari\'et\'es Riemanniennes}, J.~Lafontaine and P.~Pansu,
editors, Cedic/Fernand Nathan, 1981.

\bibitem [Gro3] {misha-hypgroups} M.~Gromov, {\it Hyperbolic groups},
in {\it Essays in Group Theory}, S.~Gersten, editor, 75--263,
Mathematical Sciences Research Institute Publications {\bf 8},
Springer-Verlag, 1987.

\bibitem [Gro4] {misha-asyminvgroups} M.~Gromov, {\it Asymptotic
Invariants of Infinite Groups}, Volume 2 of {\it Geometric Group
Theory}, G.~Niblo and M.~Roller, editors, Cambridge University Press,
1993.

\bibitem [Gro5] {misha-cc} M.~Gromov, {\it Carnot-Carath\'eodory
spaces seen from within}, in {\it Sub-Riemannian Geometry},
A.~Bellaiche and J.-J.~Risler, editors, 79--323, Birkh\"auser, 1996.

\bibitem [Gro+] {misha-translation} M.~Gromov et al., {\it Metric
Structures for Riemannian and Non-Riemannian Spaces}, Birkh\"auser,
1999.

\bibitem [GroP] {misha-pierre} M.~Gromov and P.~Pansu, {\it Rigidity
of lattices: An introduction}, in {\it Geometric Topology: Recent
Developments (Montecatini Terme, 1990)}, P.~De Bartolomeis and
F.~Tricerri, editors, 39--137, Lecture Notes in Mathematics {\bf
1504}, Springer-Verlag, 1991.

\bibitem [Gou] {Gouvea} F.~Gouv\^ea, {\it $p$-Adic Numbers: An
Introduction}, Springer-Verlag, 1993.

\bibitem [Haj] {hajlasz} P.~Haj{\l}asz, {\it Sobolev spaces on an
arbitrary metric space}, Potential Analysis {\bf 5} (1996), 403--415.

\bibitem [HajK1] {hajlasz-koskela1} P.~Haj{\l}asz and P.~Koskela, {\it
Sobolev meets Poincar\'e}, Comptes Rendus de l'Acad\'emie des Sciences
Paris S\'er.~I Math.\ {\bf 320} (1995), 1211--1215.

\bibitem [HajK2] {hajlasz-koskela2} P.~Haj{\l}asz and P.~Koskela, {\it
Sobolev Met Poincar\'e}, Memoirs of the American Mathematical Society
{\bf 688}, 2000.

\bibitem [HanH] {hanson-heinonen} B.~Hanson and J.~Heinonen, {\it An
$n$-dimensional space that admits a Poincar\'e inequality but has no
manifold points}, Proceedings of the American Mathematical Society
{\bf 128} (2000), 3379--3390.

\bibitem [Har] {Harpe} P.~de la Harpe, {\it Topics in Geometric Group
Theory}, University of Chicago Press, 2000.

\bibitem [Heb1] {Hebey1} E.~Hebey, {\it Sobolev Spaces on Riemannian
Manifolds}, Lecture Notes in Mathematics {\bf 1635}, Springer-Verlag,
1996.

\bibitem [Heb2] {Hebey2} E.~Hebey, {\it Nonlinear Analysis on
Manifolds: Sobolev Spaces and Inequalities}, Courant Lectuure Notes in
Mathematics {\bf 5}, American Mathematical Society and Courant
Institute of Mathematical Sciences, 1999.

\bibitem [Hei1] {Heinonen1} J.~Heinonen, {\it Calculus on Carnot
groups}, in {\it Fall School in Analysis (Jyv\"askyl\"a, 1994)},
1--31, Reports of the Department of Mathematics, University of
Jyv\"askyl\"a {\bf 68}, 1995.

\bibitem [Hei2] {Heinonen2} J.~Heinonen, {\it Lectures on Analysis on
Metric Spaces}, Springer-Verlag, 2001.

\bibitem [HeiK1] {HK1} J.~Heinonen and P.~Koskela, {\it Definitions of 
quasiconformality}, Inventiones Mathematicae {\bf 120} (1995), 61--79.

\bibitem [HeiK2] {HK2} J.~Heinonen and P.~Koskela, {\it From local to
global in quasiconformal structures}, Proceedings of the National
Academy of Sciences (U.S.A.) {\bf 93} (1996), 554--556.

\bibitem [HeiK3] {HK3} J.~Heinonen and P.~Koskela, {\it Quasiconformal
maps in metric spaces with controlled geometry}, Acta Mathematica {\bf
181} (1998), 1--61.

\bibitem [HeiKST] {HKST} J.~Heinonen, P.~Koskela, N.~Shanmugalingam,
and J.~Tyson, {\it Sobolev classes of Banach space-valued functions
and quasiconformal mappings}, Journal d'Analyse Math\'ematique, 
{\bf 85} (2001), 87--139.

\bibitem [HeiS] {HS} J.~Heinonen and S.~Semmes, {\it Thirty-three yes
or no questions about mappings, measures, and metrics}, Conformal
Geometry and Dynamics (an electronic journal of the American
Mathematical Society) {\bf 1} (1997), 1--12.

\bibitem [KalS] {KaS} S.~Kallunki and N.~Shanmugalingam, {\it Modulus
and continuous capacity}, Annales Academi{\ae} Scientiarum Fennic{\ae}
Mathematica Ser.\ A.I.\ Math.\ {\bf 26} (2001), 455--464.

\bibitem [Kir] {Kirchheim} B.~Kirchheim, {\it Rectifiable metric
spaces: Local structure and regularity of the Hausdorff measure},
Proceedings of the American Mathematical Society {\bf 121} (1994),
113--123.

\bibitem [KirM] {Kirchheim-Magnani} B.~Kirchheim and V.~Magnani, {\it
A counterexample to metric differentiability}, to appear, Proceedings
of the Edinburgh Mathematical Society.

\bibitem [Kig] {Kigami} J.~Kigami, {\it Analysis on Fractals},
Cambridge University Press, 2001.

\bibitem [KinS] {KinS} J.~Kinnunen and N.~Shanmugalingam, {\it
Regularity of quasi-minimizers on metric spaces}, Manuscripta
Mathematica {\bf 105} (2001), 401--423.

\bibitem [KorR] {Koranyi-Reimann} A.~Kor\'anyi and H.~Reimann, {\it
Quasiconformal mappings on the Heisenberg group}, Inventiones
Mathematicae {\bf 80} (1985), 309--338.

\bibitem [Kos] {Koskela} P.~Koskela, {\it Upper gradients and
Poincar\'e inequalities}, in {\it Lecture Notes on Analysis on Metric
Spaces}, L.~Ambrosio and F.~Serra Cassano, editors, 55--69, Scuola
Normale Superiore, Pisa, 2000.

\bibitem [KosST] {KoST} P.~Koskela, N.~Shanmugalingam, and H.~Tuominen,
{\it Removable sets for the Poincar\'e inequality on metric spaces},
Indiana University Mathematics Journal {\bf 49} (2000), 333--352.

\bibitem [Laa1] {Laakso1} T.~Laakso, {\it Ahlfors $Q$-regular spaces
with arbitrary $Q > 1$ admitting weak Poincar\'e inequality},
Geometric and Functional Analysis {\bf 10} (2000), 111--123.

\bibitem [Laa2] {Laakso2} T.~Laakso, {\it Lipschitz mappings with no
bilipschitz tangents in spaces with decent calculus}, Reports of the
Department of Mathematics, University of Helsinki {\bf 213}, 1999.

\bibitem [Laa3] {Laakso3} T.~Laakso, {\it Plane with $A_\infty$-weighted
metric not bilipschitz embeddable to ${\bf R}^n$}, Reports of the
Department of Mathematics, University of Helsinki {\bf 258}, 2000.

\bibitem [Laa4] {Laakso4} T.~Laakso, {\it Look-down equivalence
without BPI equivalence}, Reports of the Department of Mathematics,
University of Helsinki {\bf 320}, 2002.

\bibitem [LeoR] {LR} G.~Leonardi and S.~Rigot, {\it Isoperimetric
sets on Carnot groups}, to appear, Houston Journal of Mathematics.

\bibitem [LewP] {LP} H.~Lewis and C.~Papadimitriou, {\it Elements
of the Theory of Computation}, second edition, Prentice Hall, 1998.

\bibitem [LutYZ] {LYZ} E.~Lutwak, D.~Yang, and G.~Zhang, {\it Sharp
affine $L_p$ Sobolev inequalities}, preprint.

\bibitem [Luu] {Luukkainen} J.~Luukkainen, {\it Assouad dimension:
Antifractal metrization, porous sets, and homogeneous measures},
Journal of the Korean Mathematical Society {\bf 35} (1998), 23--76.

\bibitem [Mag] {Magnani} V.~Magnani, {\it Differentiability and area
formula on stratified Lie groups}, Houston Journal of Mathematics
{\bf 27} (2001), 297--323.

\bibitem [Mir] {Miranda} M.~Miranda, {\it Functions of bounded
variation on ``good'' metric measure spaces}, preprint.

\bibitem [Mon] {Monti} R.~Monti, {\it Distances, boundaries and
surface measures in Carnot--Carath\'eodory spaces}, Ph.D. thesis,
University of Trento, 2001.

\bibitem [MonM] {Monti-Morbidelli} R.~Monti and D.~Morbidelli, {\it
Regular domains in homogeneous groups}, preprint, 2001.

\bibitem [Mos] {Mostow} G.~Mostow, {\it Strong Rigidity of Locally
Symmetric Spaces}, Annals of Mathematical Studies {\bf 78}, Princeton
University Press, 1973.

\bibitem [Pan1] {Pansu1} P.~Pansu, {\it Une in\'egalit\'e
isop\'erim\'etrique sur le groupe de Heisenberg}, Comptes Rendus de
l'Acad\'emie des Sciences Paris S\'er.~I Math.\ {\bf 295} (1982),
127--130.

\bibitem [Pan2] {Pansu2} P.~Pansu, {\it An isoperimetric inequality on
the Heisenberg group}, Conference on differential geometry on
homogeneous spaces, Rend. Sem. Mat. Univ. Politec. Torino Special
Issue (1983), 159--174.

\bibitem [Pan3] {Pansu3} P.~Pansu, {\it Dimension conforme et sph\`ere
\`a l'infini des vari\'et\'es \`a courbure n\'egative}, Annales
Academi{\ae} Scientiarum Fennic{\ae} Mathematica Ser.\ A.I.\ Math.\
{\bf 14} (1989), 177--212.

\bibitem [Pan4] {Pansu4} P.~Pansu, {\it M\'etriques de
Carnot--Carath\'eodory et quasiisom\'etries des espaces sym\'etriques
de rang un}, Annals of Mathematics (2) {\bf 129} (1989), 1--60.

\bibitem [Pap] {P} C.~Papadimitriou, {\it Computational Complexity},
Addison-Wesley, 1994.

\bibitem [Sal] {Saloff-Coste} L.~Saloff-Coste, {\it A note on
Poincar\'e, Sobolev, and Harnack inequalities}, International
Mathematics Research Notices No.~2 (1992), 27--38.

\bibitem [Sem1] {Se-ind1} S.~Semmes, {\it Differentiable function
theory on hypersurfaces in ${\bf R}^n$ (without bounds on their
smoothness)}, Indiana University Mathematics Journal {\bf 39} (1990),
985-1004.

\bibitem [Sem2] {Se-finn} S.~Semmes, {\it Bilipschitz mappings and
strong $A_\infty$ weights}, Annales Academi{\ae} Scientiarum
Fennic{\ae} Mathematica Ser.\ A.I.\ Math.\ {\bf 18} (1993), 211-248.

\bibitem [Sem3] {Se-qs} S.~Semmes, {\it Good metric spaces without
good parameterizations}, Revista Matem\'atica Iberoamericana {\bf 12}
(1996), 187-275.

\bibitem [Sem4] {Se-bil} S.~Semmes, {\it On the nonexistence of
bilipschitz parameterizations and geometric problems about $A_\infty$
weights}, Revista Matem\'atica Iberoamericana {\bf 12} (1996),
337-410.

\bibitem [Sem5] {Se-top} S.~Semmes, {\it Finding Curves on General
Spaces through Quantitative Topology with Applications for Sobolev and
Poincar\`e inequalities}, Selecta Mathematica (N.S.) {\bf 2} (1996),
155-295.

\bibitem [Sem6] {Se-barc} S.~Semmes, {\it Some remarks about about
metric spaces, spherical mappings, functions and their derivatives},
Publicacions Matem\`atiques {\bf 40} (1996), 411-430.

\bibitem [Sem7] {Se-app} S.~Semmes, {\it Metric spaces and mappings
seen at many scales}, appendix in \cite{misha-translation}.

\bibitem [Sem8] {Se-cal} S.~Semmes, {\it Analysis on metric spaces},
in {\it Harmonic Analysis and Partial Differential Equations: Essays
in Honor of Alberto Calder\'on}, M.~Christ, C.~Kenig, and C.~Sadosky,
editors, 285--294, University of Chicago Press, Chicago, 1999.

\bibitem [Sem9] {Se-Trento} S.~Semmes, {\it Derivatives and difference
quotients for Lipschitz or Sobolev functions on various spaces}, in
{\it Lecture Notes on Analysis on Metric Spaces}, L.~Ambrosio and
F.~Serra Cassano, editors, 71--103, Scuola Normale Superiore, Pisa,
2000.

\bibitem [Sem10] {Se-frac} S.~Semmes, {\it Some Novel Types of Fractal
Geometry}, Oxford University Press, 2001.

\bibitem [Sha1] {Shanmugalingam1} N.~Shanmugalingam, {\it Newtonian
spaces: An extension of Sobolev spaces to metric measure spaces},
Revista Matem\'atica Iberoamericana {\bf 16} (2000), 243--279.

\bibitem [Sha2] {Shanmugalingam2} N.~Shanmugalingam, {\it Harmonic
functions on metric spaces}, Illinois Journal of Mathematics {\bf 45}
(2001), 1021--1050.

\bibitem [Ste1] {Stein1} E.~Stein, {\it Singular Integrals and
Differentiability Properties of Functions}, Princeton University
Press, 1970.

\bibitem [Ste2] {Stein2} E.~Stein, {\it Harmonic Analysis:
Real-Variable Methods, Orthogonality, and Oscillatory Integrals},
Princeton University Press, 1993.

\bibitem [SteW] {SW} E.~Stein and G.~Weiss, {\it Introduction to
Fourier Analysis on Euclidean Spaces}, Princeton University Press,
1971.

\bibitem [Tho] {Thompson} A.~Thompson, {\it Minkowski Geometry},
Cambridge University Press, 1996.

\bibitem [Tys1] {Tyson1} J.~Tyson, {\it Quasiconformality and
quasisymmetry in metric measure spaces}, Annales Academi{\ae}
Scientiarum Fennic{\ae} Mathematica Ser.\ A.I.\ Math.\ {\bf 23}
(1998), 525--548.

\bibitem [Tys2] {Tyson2} J.~Tyson, {\it Metric and geometric
quasiconformality in Ahlfors regular Loewner spaces}, Conformal
Geometry and Dynamics (an electronic journal of the American
Mathematical Society) {\bf 5} (2001), 21--73.

\bibitem [V\"ai] {Vaisala} J.~V\"ais\"al\"a, {\it Metric duality in
Euclidean spaces}, Mathematica Scandinavica {\bf 80} (1997), 249--288.

\bibitem [VarSC] {VSC} N.~Varopoulos, L.~Saloff-Coste, and T.~Coulhon,
{\it Analysis and Geometry on Groups}, Cambridge University Press, 1992.

\bibitem [Whe] {Wheeden} R.~Wheeden, {\it Some weighted Poincar\'e
estimates in spaces of homogeneous type}, in {\it Lecture Notes on
Analysis on Metric Spaces}, L.~Ambrosio and F.~Serra Cassano, editors,
105--121, Scuola Normale Superiore, Pisa, 2000.


\end{thebibliography}
\end{document}